\documentclass[sn-mathphys-num]{sn-jnl}

\usepackage{graphicx}%
\usepackage{multirow}%
\usepackage{amsmath,amssymb,amsfonts}%
\usepackage{amsthm}%
\usepackage{mathrsfs}%
\usepackage[title]{appendix}%
\usepackage{xcolor}%
\usepackage{textcomp}%
\usepackage{manyfoot}%
\usepackage{booktabs}%
\usepackage{algorithm}%
\usepackage{algorithmicx}%
\usepackage{algpseudocode}%
\usepackage{listings}%

\theoremstyle{thmstyleone}%

\newtheorem{lemma}{Lemma}

\newtheorem{theorem}{Theorem}
\numberwithin{equation}{section}

\theoremstyle{thmstyletwo}%

\theoremstyle{thmstylethree}%

\begin{document}

\title[Wright Operators]{Some properties of Wright Operators}


\author*{\fnm{Prashantkumar G. } \sur{Patel}}\email{prashant225@gmail.com;prashant225@spuvvn.edu}

%

\affil*{\orgdiv{Department of Mathematics}, \orgname{Sardar Patel University}, \orgaddress{\street{}, \city{Vallabh Vidyanagar}, \postcode{388 120}, \state{Gujarat}, \country{India}}
orchid ID: \url{https://orcid.org/0000-0002-8184-1199 }}




\abstract{This note aims to present a new sequence of positive linear operators involving the Wright function. Furthermore, the present research established the moments of these newly defined operators and estimated the convergence rate using the classical modulus of continuity. Additionally, the convergence rate in the Lipschitz spaces, their $A$-statistical convergence property have been covered.}




\maketitle

\section{Introduction}\label{sec1}
 A new generalization of Bernstein operators was introduced periodically in approximation theory (see \cite{totik1994approximation,lewanowicz2004generalized,sinha2009inverse,gupta2013approximation,gupta2021operators,mishra2014generalized}). Many positive linear operators are derived from special functions. For instance, the operators were introduced and studied in \cite{ozarslan2013astatistical} using the Mittag-Leffler function, the operators were studied in \cite{gupta2020convergence} with the help of gamma function, and using Laguerre polynomials, the positive linear operators introduced by Gupta \cite{gupta2024new}. This note aims to merge two mathematical topics-approximation theory and special functions. \\
The special function  
\begin{equation}\label{wrightfunction}
	\phi(\rho,\beta;z)=\phi_{\rho,\beta}(z) = \sum_{k=0}^{\infty} \frac{z^k}{k!\Gamma(\rho k +\beta)}, ~~~(\rho>-1,\beta,z\in \mathbb{C})
\end{equation}
named after the British mathematician E. M. Wright, has appeared for the first time in the case $\rho>0$ in connection with his investigations in the asymptotic theory of partitions in \cite{wright1933coefficients}. \\
Note that $\phi(1,1;z)= I_0(2\sqrt{z})$,\\
 where  $I_n(z)=i^{-n}J_n(iz)=\dfrac{(z/2)^n}{\Gamma(1+n)} {}_{0}F_{1}\left(-;1+n;\dfrac{z^2}{4}\right)$, $n$ not a negative integer. The function $I_n(z)$ is called the modified Bessel function of the first kind of index $n$. Here, $J_n$ is Bessel function and ${}_{0}F_{1}$ is hypergeometric function. The modified function $I_n$ is related to $J_n$ in much the same way that the hyperbolic function is related to the trigonometric function. \\
Also, note that 
$$\phi_{1,2}(z)=\dfrac{1}{\sqrt{z}}I_1\left(2\sqrt{z}\right);~~~~~ \phi_{1,m}(z)=\left(\dfrac{1}{\sqrt{z}}\right)^{m/2} I_{m-1}\left(2\sqrt{z}\right).$$
It was demonstrated in \cite{gray1895treatise,bateman1953higher} that, the modified Bessel function $I_n(z)$ of the first kind does not have zeros in the right half-plane. \cite{khanmamedov2010zeros} contains an analysis of $n$-zeros of the function $I_n(z)$ in the left half-plane. Since $I_n(z)$ has no zeros in $[0,\infty)$, the Wright function $\phi_{1,m}(z)$ also has no zeros in $[0,\infty)$.\\
Let $\beta>1$ be fixed. We provide a novel class of positive linear operators involving the Wright function for every $n\in \mathbb{N}$ as 
\begin{equation} \label{wrightoperators}
	W_n^{(\beta)}(f;x)=\dfrac{1}{\phi_{1,\beta}\left(nx\right)}\sum_{k=0}^{\infty} f\left(\dfrac{k+\beta}{n}\right)\dfrac{(nx)^k}{k! \Gamma(k+\beta)}, 
\end{equation}
where $	f\in E :=\left\{f\in C\left[0,\infty\right): \lim_{x\to \infty} \dfrac{f(x)}{1+x^2} \textrm{ is finite }\right\}$ and $C\left[0,\infty\right)$ denote the space of continuous functions defined on $[0,\infty)$. Recall that the Banach lattice $E$ is endowed with the norm 
$$\|f\|_2 :=\sup_{x\in [0,\infty)} \dfrac{|f(x)|}{1+x^2}.$$  
One can note that the operators $W_n^{(\beta)}$ defined in \eqref{wrightoperators} are linear and positive. We called it the Wright operators. \\
\indent Research on approximating continuous signals using a sequence of positive linear operators is ongoing. We provide other findings along the same lines in this note, but with the aid of recently established positive linear operators $W_n^{(\beta)}$. There is no prior literature that introduced these operators \eqref{wrightoperators}. We estimate the moments and central moments of the operators \eqref{wrightoperators} up to the fourth order in the second section. The rate of convergence of these positive linear operators is discussed in the next section. In the same section 3, we have proved that $W_n^{(\beta)}$ maps $E$ into itself. We have also determined the rate of convergence for these operators by utilizing the modulus of continuity. For the operators $W_n^{(\beta)}$, we have established a statistical Voronovskaya-type theorem in section 4.
\section{Some Lemmas}
To analyze approximation features of the operators \eqref{wrightoperators}, a few inequalities for the Wright function are required.
\begin{lemma}\cite[Theorem 6.1]{mehrez2017functional}\label{wrightfunprop}
	Let $\alpha,\beta>0$. Then the following assertions are true:
	\begin{equation}
		\Gamma(\beta+\alpha)\phi_{\alpha,\beta+\alpha} (z)\leq \Gamma(\beta)\phi_{\alpha,\beta}(z), \textrm{ for any }z>0.
	\end{equation}
\end{lemma}
Direct calculations allow one to declare the following lemma:
\begin{lemma}\label{moments1} Let $\phi_x^2\left(t\right)=\left(t-x\right)^2$, for each $x\geq 0$, $\beta>1$ and $n \in \mathbb{N}$, we have
	\begin{enumerate} 
		\item $	W_n^{(\beta)}(1;x)=1;$
		\item $\left|W_n^{(\beta)}(t;x)-x\right|\leq \dfrac{\beta}{n};$
		\item $\left|W_n^{(\beta)}(t^2;x)-x^2\right|\leq  \dfrac{x(1+2\beta)}{\beta n}+\dfrac{\beta^2}{n^2}.;$
		\item $\left|W_n^{(\beta)}(t^3;x)-x^3\right| \leq \dfrac{3x^2}{n\beta}+\dfrac{x(1+3\beta+\beta^2)}{\beta n^2}+
		\dfrac{\beta^3}{n^3};$
		\item $\left|W_n^{(\beta)}(t^4;x)-x^4\right| \leq \frac{(4\beta +6)x^3}{n\beta(\beta+1)(\beta+2)} +\frac{(6\beta^2+12\beta+7)x^2}{n^2\beta(\beta+1)}+\frac{(4\beta^3+6\beta^2+4\beta+1)x}{n^3\beta}+\frac{\beta^4}{n^4}.$
	\end{enumerate}
\end{lemma}
\begin{proof}
Since, 
$$\sum_{k=0}^{\infty} \frac{(nx)^k}{k!\Gamma(k+\beta)}=\phi_{1,\beta}(nx),$$
one can have $W_n^{(\beta)}(1,x)=1$. Using Lemma \ref{wrightfunprop}, $\Gamma(\beta+1)\phi_{1,\beta+1}\left(nx\right)\leq\Gamma(\beta)\phi_{1,\beta}\left(nx\right)$ for any $x\in [0,\infty)$ and $n\in \mathbb{N}$, we get 
\begin{align*}
	W_n^{(\beta)}(t;x)&=\dfrac{1}{\phi_{1,\beta}\left(nx\right)}\sum_{k=0}^{\infty} \dfrac{1}{n}\dfrac{(nx)^k}{(k-1)! \Gamma(k+\beta)}+\dfrac{\beta}{n}\\
	&= \dfrac{x}{\phi_{1,\beta}\left(nx\right)}\sum_{k=0}^{\infty} \dfrac{(nx)^{k}}{k! \Gamma(k+\beta+1)}+\dfrac{\beta}{n}=x\dfrac{\phi_{1,\beta+1}(nx)}{\phi_{1,\beta}(nx)}+\dfrac{\beta}{n}\leq  x+\dfrac{\beta}{n} \textrm{ for } \beta>1.
\end{align*}
Again applying Lemma \ref{wrightfunprop}, $\Gamma(\beta+2)\phi_{1,\beta+2}\left(nx\right)\leq\Gamma(\beta+1)\phi_{1,\beta+1}\left(nx\right)\leq \Gamma(\beta)\phi_{1,\beta}\left(nx\right)$ for any $x\in [0,\infty)$ and $n\in \mathbb{N}$, we get 	
\begin{align*}
	W_n^{(\beta)}(t^2;x)&=\dfrac{1}{\phi_{1,\beta}\left(nx\right)}\sum_{k=1}^{\infty} \dfrac{(k+\beta)^2}{n^2}\dfrac{(nx)^k}{(k-1)! \Gamma(k+\beta)}\\
	&= \frac{1}{n^2 \, \phi_{1,\beta}(nx)} \Bigg[ 
	(nx)^2 \sum_{k=0}^{\infty} \frac{(nx)^k}{k! \, \Gamma(k+2+\beta)} 
	+ (1+2\beta) \, nx \sum_{m=0}^{\infty} \frac{(nx)^m}{m! \, \Gamma(m+1+\beta)} 
	+ \beta^2 \, \phi_{1,\beta}(nx) 
	\Bigg] \nonumber \\
	&= \frac{1}{n^2 \, \phi_{1,\beta}(nx)} \Big[(nx)^2 \, \phi_{1,\beta+2}(nx) + (1+2\beta) \, nx \, \phi_{1,\beta+1}(nx) + \beta^2 \, \phi_{1,\beta}(nx) \Big].
\end{align*}
Since $\beta>1$, $$\left|	W_n^{(\beta)}(t^2;x)-x^2 \right|\leq \dfrac{x(1+2\beta)}{\beta n}+\dfrac{\beta^2}{n^2}.$$
Now, 
\begin{eqnarray*}
	W_n^{(\beta)}(t^3;x)
		&=& \frac{1}{n^3\phi_{1,\beta}(nx)}\sum_{k=0}^{\infty}(k+\beta)^3\frac{(nx)^k}{k!,\Gamma(k+\beta)}\\
		&=& \frac{1}{n^3\phi_{1,\beta}(nx)}\Big[ (nx)^3\phi_{1,\beta+3}(nx)
		+3(1+\beta)(nx)^2\phi_{1,\beta+2}(nx)\\
		&&+(1+3\beta+3\beta^2)nx\phi_{1,\beta+1}(nx)
		+\beta^3\phi_{1,\beta}(nx)\Big].\\
		&=&x^3\dfrac{\phi_{1,\beta+3}\left(nx\right)}{\phi_{1,\beta}\left(nx\right)}+3(1+\beta)x^2 \dfrac{\phi_{1,\beta+2}\left(nx\right)}{n\phi_{1,\beta}\left(nx\right)}+x(1+3\beta+\beta^2)\dfrac{\phi_{1,\beta+1}\left(nx\right)}{n^2\phi_{1,\beta}\left(nx\right)}+\dfrac{\beta^3}{n^3}\\
	&\leq&  \left|\dfrac{x^3\Gamma(\beta)}{\Gamma(\beta+3)}\right| + \left|\dfrac{3(1+\beta)x^2\Gamma(\beta)}{n\Gamma(\beta+2)} \right|+\left|\dfrac{x(1+3\beta+\beta^2)\Gamma(\beta)}{n^2\Gamma(\beta+1)}\right|\\
	&=& \dfrac{x^3}{\beta(\beta+1)(\beta+2)} +\dfrac{3x^2}{n\beta}+\dfrac{x(1+3\beta+\beta^2)}{\beta n^2}+
	\dfrac{\beta^3}{n^3}. 
\end{eqnarray*}
Since, $\beta>1$, 
$$\displaystyle \left|W_n^{(\beta)}(t^3;x)-x^3\right| \leq \dfrac{3x^2}{n\beta}+\dfrac{x(1+3\beta+\beta^2)}{\beta n^2}+
\dfrac{\beta^3}{n^3}.$$
Similarly,
\begin{eqnarray*}
	W_n^{(\beta)}(t^4;x)&=&\dfrac{1}{\phi_{1,\beta}\left(nx\right)}\sum_{k=0}^{\infty} \dfrac{(k+\beta)^4}{n^4}\dfrac{(nx)^k}{k! \Gamma(k+\beta)}\\
	&=& \frac{1}{n^4\phi_{1,\beta}(nx)}\Big[ (nx)^4\phi_{1,\beta+4}(nx)+(4\beta+6)(nx)^3\phi_{1,\beta+3}(nx)\\
	&& +(6\beta^2+12\beta+7)(nx)^2\phi_{1,\beta+2}(nx)+(4\beta^3+6\beta^2+4\beta+1)nx\phi_{1,\beta+1}(nx) +\beta^4\phi_{1,\beta}(nx)\Big].\\
	&=&x^4\dfrac{\phi_{1,\beta+4}\left(nx\right)}{\phi_{1,\beta}\left(nx\right)}+(4\beta+6)x^3\dfrac{\phi_{1,\beta+3}\left(nx\right)}{n\phi_{1,\beta}\left(nx\right)}+(6\beta^2+12\beta+7)x^2\dfrac{\phi_{1,\beta+2}\left(nx\right)}{n^2\phi_{1,\beta}\left(nx\right)}\\
	&&+(4\beta^3+6\beta^2+4\beta+1)x\dfrac{\phi_{1,\beta+1}\left(nx\right)}{n^3\phi_{1,\beta}\left(nx\right)}+\dfrac{\beta^4}{n^4}\\
	&\leq& \left|\dfrac{ x^4 \Gamma(\beta)}{\Gamma(\beta+4)}\right|+\left| \dfrac{(4\beta+6)x^3\Gamma(\beta)}{n\Gamma(\beta+3)}\right|+\left|\dfrac{(6\beta^2+12\beta+7)x^2\Gamma(\beta)}{n^2\Gamma(\beta+2)}\right|\\
	&&+\left|\dfrac{(4\beta^3+6\beta^2+4\beta+1)x\Gamma(\beta)}{n^3\Gamma(\beta+1)}\right|+\dfrac{\beta^4}{n^4}\\
	&=&  \dfrac{ x^4}{\beta(\beta+1)(\beta+2)(\beta+3)}+\dfrac{(4\beta +6)x^3}{n\beta(\beta+1)(\beta+2)} +\dfrac{(6\beta^2+12\beta+7)x^2}{n^2\beta(\beta+1)}\\
	&&+\dfrac{(4\beta^3+6\beta^2+4\beta+1)x}{n^3\beta}+\dfrac{\beta^4}{n^4}.
\end{eqnarray*}
Since, $\beta>1$, 
$$\left|W_n^{(\beta)}(t^4;x)-x^4\right| \leq \dfrac{(4\beta +6)x^3}{n\beta(\beta+1)(\beta+2)} +\dfrac{(6\beta^2+12\beta+7)x^2}{n^2\beta(\beta+1)}+\dfrac{(4\beta^3+6\beta^2+4\beta+1)x}{n^3\beta}+\dfrac{\beta^4}{n^4}.$$
\end{proof}
\begin{lemma}
	Let $\phi_x^i\left(t\right)=\left(t-x\right)^i$, $i\in \mathbb{N}$, for each $x\geq 0$, $\beta>1$ and $n \in \mathbb{N}$, we have
	\begin{enumerate} 
		\item $W_n^{(\beta)}(\phi_x^1;x)\leq \dfrac{\beta}{n}$
		\item $W_n^{(\beta)}(\phi_x^2;x)\leq \frac{x(1+2\beta)}{\beta n} + \frac{2x\beta}{n} + \frac{\beta^2}{n^2}.$
		\item $W_n^{(\beta)}(\phi_x^3;x) \leq \frac{3x^2}{n\beta}+\frac{x(1+3\beta+\beta^2)}{\beta n^2}+\frac{\beta^3}{n^3} 
		+ 3x\left(\frac{x(1+2\beta)}{\beta n}+\frac{\beta^2}{n^2}\right) 
		+ \frac{3x^2\beta}{n} .$
		\item $W_n^{(\beta)}(\phi_x^4;x) \leq \left( \frac{4\beta+6}{\beta(\beta+1)(\beta+2)} + \frac{18}{\beta} + 12 + 4\beta \right)\frac{x^3}{n} 
		+ \left( \frac{6\beta^2+12\beta+7}{\beta(\beta+1)} + \frac{4}{\beta} + 12 + 4\beta + 6\beta^2 \right)\frac{x^2}{n^2} 
		+ \left( 4\beta^3 + 4\beta^2 + 6\beta + 4 + \frac{1}{\beta} \right)\frac{x}{n^3} 
		+ \frac{\beta^4}{n^4}.$
	\end{enumerate}
\end{lemma}
\begin{proof}
Using Lemma \ref{moments1} and simple computation gives $W_n^{(\beta)}(\phi_x^1;x)=0$ and 
$$W_n^{(\beta)}(\phi_x^2;x)\leq \left|W_n^{(\beta)}(t^2;x)-x^2\right|+2x\left|W_n^{(\beta)}(t;x)-x\right|+x^2\left|W_n^{(\beta)}(1;x)-1\right|\leq \frac{x(1+2\beta)}{\beta n} + \frac{2x\beta}{n} + \frac{\beta^2}{n^2}.$$
Further, 
\begin{align*}
	W_n^{(\beta)}(\phi_x^3;x)&\leq \left|W_n^{(\beta)}(t^3;x)-x^3\right|+3x\left|W_n^{(\beta)}(t^2;x)-x^2\right|+3x^2\left|W_n^{(\beta)}(t;x)-x\right|+x^3 \left|W_n^{(\beta)}(1;x)-1\right|\\
	&\leq \dfrac{3x^2}{n\beta(\beta+1)}+\dfrac{x}{n^2\beta}+3x\left( \dfrac{x}{n\beta }\right)=\dfrac{3 x^2(\beta +2)}{n\beta(\beta +1) }+\dfrac{x}{n^2\beta}.
\end{align*}
Finally, we have 
\begin{align*}
	W_n^{(\beta)}(\phi_x^4;x)&\leq \left|W_n^{(\beta)}(t^4;x)-x^4\right|+4x\left|W_n^{(\beta)}(t^3;x)-x^3\right|\\
	&+6x^2\left|W_n^{(\beta)}(t^2;x)-x^2\right|+4x^3 \left|W_n^{(\beta)}(t;x)-x\right|+x^4\left|W_n^{(\beta)}(1;x)-1\right|\\
	&\leq \left(\frac{4\beta+6}{\beta(\beta+1)(\beta+2)}+\frac{12}{\beta}+\frac{6(1+2\beta)}{\beta}+4\beta\right)\frac{x^3}{n} \\
	&
	+ \left(\frac{6\beta^2+12\beta+7}{\beta(\beta+1)} + \frac{4(1+3\beta+\beta^2)}{\beta} + 6\beta^2\right)\frac{x^2}{n^2} 
	+ \left(\frac{4\beta^3+6\beta^2+4\beta+1}{\beta}+4\beta^3\right)\frac{x}{n^3} 
	+ \frac{\beta^4}{n^4}.
\end{align*}
\end{proof}
\section{Rate of Convergence}
In this section, the rate of convergence of $W_n^{(\beta)}$ is discussed. The following lemma proves that $W_n^{(\beta)}$ maps $E$ into itself.
\begin{lemma}
	Let $\beta>1$ be fixed, then there exists a constant $M(\beta)$ such that,  
	$$\omega(x) W_n^{(\beta)} \left(\dfrac{1}{\omega}, x\right) \leq M(\beta)$$
	holds for all $x\in [0,\infty)$, $n\in \mathbb{N}$ and  $\omega(x)=\dfrac{1}{1+x^2}$. Furthermore, for all $f\in E$, we have
	$$\|D_n^{(\beta)}(f, \cdot)\|_2\leq M(\beta)\|f\|_2$$
\end{lemma}
\begin{proof} From Lemma \ref{moments1}, we have 
\begin{eqnarray*}
	\omega(x) W_n^{(\beta)}\left(\dfrac{1}{\omega},x\right) &=& \dfrac{1}{1+x^2} \left[W_n^{(\beta)}\left(1,x\right) + W_n^{(\beta)}\left(t^2,x\right)\right]\\
	&\leq & \dfrac{1}{1+x^2} \left[ 1+x^2 +\dfrac{x(1+2\beta)}{\beta n}+\dfrac{\beta^2}{n^2}\right]\\
	&\leq& M(\beta)
\end{eqnarray*}
This follows, following inequality
$$\omega(x) \left|W_n^{(\beta)}\left(f,x\right)\right|=\omega(x) \left|W_n^{(\beta)}\left(\omega \dfrac{f}{\omega},x\right)\right|\leq \|f\|_2 \omega(x)W_n^{(\beta)}\left( \dfrac{1}{\omega},x\right)\leq M(\beta) \|f\|_2 $$
Taking supremum over $x\in [0,\infty)$ in the above inequality, gives the result.
\end{proof}
Now, recall that the usual modulus of continuity of $f$ on the closed interval $[0,B]$ is defined by
$$\omega_B(f,\delta) =\sup\{ |f(t)-f(x)|: |t-x|\leq \delta, ~~x,t\in [0,B]\}$$
It is well known that, for a function $f \in E$, we have $\displaystyle \lim_{\delta\to \infty}	\omega_B(f,\delta)=0$.
The next theorem gives the rate of convergence of the operators $W_n^{(\beta)}(f,x)$, for all $f \in E$.\\
\begin{theorem}\label{theorem1}
	Let $\beta>1$, $f\in E$ and $\omega_{B+1}(f,\delta)$, $(B>0)$ be its modulus of continuity on the finite interval $[0,B+1]\subset [0,\infty)$, then 
	$$\|W_n^{(\beta)}\left(f,\cdot\right)-f(\cdot)\|_{C[0,B]} \leq M_f(\beta,B) \delta_n(\beta,B)+ 2\omega_{B+1}\left(f, \delta_n^{1/2}(\beta,B)\right)$$
	where $\delta_n(\beta,B)=\dfrac{B(1+2\beta)}{\beta n}+\dfrac{\beta^2}{n^2}$	and $M_f(\beta,B)$ is an absolute constant depending on $f$, $\beta$ and $B$. 
\end{theorem}
\begin{proof} Let $\beta >0$ be fixed. For $x\in [0,B]$ and $t\leq B+1$, we have following well-known inequality 
\begin{equation}\label{eq1}
	|f(t)-f(x)| \leq \omega_{B+1}(f, |t-x|) \leq \left(1+ \dfrac{|t-x|}{\delta}\right) \omega_{B+1}(f,\delta)
\end{equation}
where $\delta>0$. Now, for $x\in [0,B]$ and $t>B+1$, using the fact that $t-x>1$, we have 
\begin{eqnarray}\label{eq21}
	|f(t)-f(x)|&\leq& A_f (1+x^2 +t^2)\notag \\
	& \leq &  A_f (2 + 3x^2 +2(t-x)^2 )\notag \\
	&\leq & 6 A_f (1+B^2) (t-x)^2 
\end{eqnarray}
Using \eqref{eq1} and \eqref{eq21}, we get for all $x\in [0,B]$ and $t\geq 0$, we get  
\begin{eqnarray}\label{eq2}
	|f(t)-f(x)|&\leq& 6 A_f (1+B^2) (t-x)^2 + \left(1+ \dfrac{|t-x|}{\delta}\right) \omega_{B+1}(f,\delta)
\end{eqnarray}
Therefore, 
\begin{eqnarray*}
	\left|W_n^{(\beta)}\left(f,x\right)-f(x)\right|\leq 6 A_f (1+B^2)W_n^{(\beta)}\left(\phi_x^2,x\right)+ \left(1+\dfrac{W_n^{(\beta)}\left(|t-x|,x\right)}{\delta} \right)\omega_{B+1}(f,\delta)
\end{eqnarray*}
Applying Cauchy-Schwarz inequality and Lemma \ref{moments1}, we get
\begin{eqnarray*}
	\left|W_n^{(\beta)}\left(f,x\right)-f(x)\right|&\leq& 6 A_f (1+B^2)W_n^{(\beta)}\left(\phi_x^2,x\right)+ \left(1+\dfrac{\left[W_n^{(\beta)}\left(\phi_x^2,x\right)\right]^{1/2}}{\delta} \right)\omega_{B+1}(f,\delta) \\
	&\leq& 	6 A_f (1+B^2)\left(\dfrac{B(1+2\beta)}{\beta n}+\dfrac{\beta^2}{n^2}\right)+ \left(1+\dfrac{\left[\dfrac{B(1+2\beta)}{\beta n}+\dfrac{\beta^2}{n^2}\right]^{1/2}}{\delta} \right)\omega_{B+1}(f,\delta)\\
	&\leq & M_f(\beta,B) \delta_n^2(\beta,B) + 2\omega_{B+1} \left(f,\left(\delta_n(\beta,B)\right)^{1/2}\right)
\end{eqnarray*}
where 
$M_f(\beta,B)=6 A_f (1+B^2)$ and $\delta_n(\beta,B)=\dfrac{B(1+2\beta)}{\beta n}+\dfrac{\beta^2}{n^2}$. Hence the proof. \end{proof}
\section{$A$-statistical Convergence}
In this section, first we discussed some definitions and notations on the concept of $A$-statistical convergence. Let $A=\left(a_{nk}\right)$, $\left(n,k\in \mathbb{N}\right)$, be a non-negative, infinite summability matrix. For a given sequence $x:=(x_k)$, the $A$-transform of $x$ denoted by $Ax: \left((Ax)_n\right)$ is defined as 
$$\left(Ax\right)_n =\sum_{k=1}^{\infty} a_{nk} x_k,$$
provided the series converges for each $n$. $A$ is said to be regular if
$\displaystyle \lim_n (Ax)_n =L$ whenever $\lim_n x_n =L$. We say that, the sequence $x=(x_n)$ is $A$-statistically convergent to $L$ and write $st_A-\lim_n x_n =L$ if for every $\epsilon>0$, $\displaystyle \lim_n \sum_{k : |x_k-L|\geq \epsilon} a_{nk}=0$.\\
Replacing $A$ by $C_1$, the Ces\'{a}ro matrix of order one, the $A$-statistical convergence reduces to statistical convergence. Similarly, if we take $A = I$, the identity matrix, then $A$-statistical convergence coincides with ordinary convergence. The statistical convergence of various types of operators has been studied by several researchers  (see \cite{mishra2017statistical,duman2008statistical,gadjiev2002some,duman2006statistical,gupta2009statistical,ispir2008statistical,gupta2014convergence,gupta2021modifications,gupta2021higher,mursaleen2009invariant,gurhan2015durrmeyer,garrancho2019general,chandra2022approximation,agrawal2022stancu,maurya2018approximation}). Now, we prove weighted Korovkin theorem via $A$-statistical convergence.\\
Here, we recall the weighted Korovkin type approximation theorem for the A-statistical convergence was given by Duman and Orhan in \cite{duman2004statistical}.
\begin{theorem}\cite{duman2004statistical}
	Let $A$ be a non-negative regular summability matrix and let $\bar{\rho}_1$; $\bar{\rho}_2$ weight functions such that $$\lim_{|x|\to \infty} \dfrac{\bar{\rho}_1(x)}{\bar{\rho}_2(x)}=0.$$
	Assume that $(T_n)_{n\geq 1}$ is a sequence of positive linear operators from $C_{\bar{\rho}_1}(\mathbb{R})$ into $B_{\bar{\rho}_2}(\mathbb{R})$, One has
	$$st_A-\lim_n \|T_nf-f\|_{\bar{\rho}_2}=0,$$
	for all $f\in C_{\bar{\rho}_1}(\mathbb{R})$ if and only if 
	$$st_A-\lim_n \|T_nF_v-F_v\|_{\bar{\rho}_2}=0, \textrm{ for all } v=0,1,2,$$
	where
	$$F_v(x) = \dfrac{x^v\bar{\rho}_1(x)}{1+x^2}, v=0,1,2.$$
\end{theorem}
By using this theorem the following Korovkin type theorem can be proved for $\left(W_n^{(\beta)}\right)$.
\begin{theorem}\label{theoremstastical}
	Let $A=(a_{nk})$ be a non-negative regular summability matrix, $\beta>1$ be fixed and $x\in [0,\infty)$, then for all $f\in E$, we have 
	$$st_A-\lim_n \|W_n^{(\beta)}\left(f,\cdot\right) -f\|_{2}=0$$
\end{theorem}
\begin{proof} 
From \cite[p. 191, Th. 3]{duman2004statistical}, it is sufficient to show that $st_A-\lim_n \|W_n^{(\beta)}\left(t^i,\cdot\right) -x^i\|_{2} =0$, where $i=0,1,2$. \\
In view of Lemma \ref{moments1}, it follows that 
$$st_A-\lim_n \|W_n^{(\beta)}\left(1,\cdot\right) -1\|_{2}=0$$
and 
$$st_A-\lim_n \|W_n^{(\beta)}\left(t,\cdot\right) -x\|_{2}=0.$$
Now, 
$$\|W_n^{(\beta)}\left(t^2,\cdot\right) -x^2\|_{2}\leq \sup_{x\geq 0}\left(\dfrac{x^2(1+2\beta)}{\beta n(1+x^2)}+\dfrac{x\beta^2}{n^2(1+x^2)}\right)\leq \dfrac{1+2\beta}{n\beta}+\dfrac{\beta^2}{n^2}$$
Given $r>0$, choose $\epsilon>0$ such that $\epsilon <r$. For fixed $\beta >1$, define the following sets: 
$$U:=\left\{n:  \|W_n^{(\beta)}\left(t,\cdot\right) -x\|_{2} \geq \epsilon\right\}$$
$$U_1:=\left\{n : \dfrac{(1+2\beta)}{\beta n}+\dfrac{\beta^2}{n^2} \geq \dfrac{\epsilon}{2}\right\}$$
Then it is clear that $U\subset U_1$, this gives 
\begin{equation}\label{theoremeq1}
	\sum_{k\in U} a_{jk} \leq \sum_{k\in U_1} a_{jk}
\end{equation}
Letting $j\to \infty$ in  \eqref{theoremeq1}, we have 	$\lim_j \sum_{k\in U} a_{jk}=0$. This proves that $st_A-\lim_n \dfrac{1}{n\beta}=0$, this also implies 
$$st_A- \lim_n \omega_{B+1} \left(f, \sqrt{\dfrac{1}{n\beta}}\right) =0$$
Using theorem \ref{theorem1}, we get desired result.\end{proof}
\begin{lemma}\label{theorem1.23}
	Let $A=(a_{in})$ be a non-negative regular summability matrix, then we have 
	$$st_A-\lim_{n\to \infty} nW_n^{(\beta)}\left((t-x)^4,x\right)=\frac{4\beta+6}{\beta(\beta+1)(\beta+2)} + \frac{18}{\beta} + 12 + 4\beta .$$
	uniformly with respect to $x\in [0, b]$ with $b >0$.
\end{lemma}
\begin{proof} Note that 
\begin{eqnarray*}
	W_n^{(\beta)}\left((t-x)^4,x\right)&\leq& \left( \frac{4\beta+6}{\beta(\beta+1)(\beta+2)} + \frac{18}{\beta} + 12 + 4\beta \right)\frac{x^3}{n} 
	+ \left( \frac{6\beta^2+12\beta+7}{\beta(\beta+1)} + \frac{4}{\beta} + 12 + 4\beta + 6\beta^2 \right)\frac{x^2}{n^2}\\
	&& + \left( 4\beta^3 + 4\beta^2 + 6\beta + 4 + \frac{1}{\beta} \right)\frac{x}{n^3} 
	+ \frac{\beta^4}{n^4}.
\end{eqnarray*}
This gives 
$$\left|nW_n^{(\beta)}\left((t-x)^4,x\right)-\left(\frac{4\beta+6}{\beta(\beta+1)(\beta+2)} + \frac{18}{\beta} + 12 + 4\beta\right) \right|\leq \dfrac{h_1(\beta)}{n}x^2+\dfrac{h_2(\beta)}{n^2}x+\dfrac{\beta^4}{n^3}$$
where $h_1$ and $h_2$ are some functions of $\beta$. \\
For $x\in [0,b]$, we have 
\begin{equation}\label{eq41}
	\left|nW_n^{(\beta)}\left((t-x)^4,x\right)-\left(\frac{4\beta+6}{\beta(\beta+1)(\beta+2)} + \frac{18}{\beta} + 12 + 4\beta\right) \right|\leq B\left\{\dfrac{1}{n}+\dfrac{1}{n^2}+\dfrac{1}{n^3}\right\}
\end{equation}
where $B=\max\{h_1(\beta)b^2, h_2(\beta)b,\beta^4\}$. Now, for a given $\epsilon > 0$, define the following sets:
\begin{eqnarray*}
	D&:=&\left\{n: \left|nW_n^{(\beta)}\left((t-x)^4,x\right)-\left(\frac{4\beta+6}{\beta(\beta+1)(\beta+2)} + \frac{18}{\beta} + 12 + 4\beta\right)\right|\geq \epsilon\right\}\\
	D_1&:=&\left\{n: \dfrac{1}{n} \geq \dfrac{\epsilon}{2B}\right\}\\
	D_2&:=&\left\{n: \dfrac{1}{n} \geq \sqrt{\dfrac{\epsilon}{2B}}\right\}\\
	D_3&:=&\left\{n: \dfrac{1}{n} \geq \sqrt[3]{\dfrac{\epsilon}{B}}\right\}
\end{eqnarray*}
Hence, by inequality \eqref{eq41}, we see that $D\subset D_1\cup D_2\cup D_3$. Then for any $j\in \mathbb{N}$, we have 
\begin{equation}
	\sum_{n\in D} a_{jn} \leq 	\sum_{n\in D_1} a_{jn}+	\sum_{n\in D_2} a_{jn} \sum_{n\in D_3} a_{jn}.
\end{equation}
Taking limit as $j \to \infty $ on the both sides of above inequality and using the fact that $\displaystyle st_A-\lim_{n\to \infty }\left(\dfrac{1}{n}\right)= 0$, we conclude that
$$\lim_{j\to \infty} \sum_{n\in D} a_{jn}=0.$$
Hence the result. \end{proof}
\begin{theorem}[Statistical Voronovskaya-type theorem for the operator $W_n^{(\beta)}$]
	Let $A=(a_{jn})$ be a nonnegative regular 	summability matrix, $\beta>1$ then for every $f\in E$ with $f',f''\in E$, we have 
	$$st_A-\lim_{n\to \infty} n\left(W_n^{(\beta)}(f,x)-f(x)\right)=\frac{1 + 2\beta + 2\beta^2}{\beta}\dfrac{x}{2}f''(x)$$
	uniformly with respect to $x\in [0,b]$ with $b>0$. 
\end{theorem}
\begin{proof}
 Let $f,f',f''\in E$ and $x\in [0,b]$. Define the function $\Phi_x$ by 
\begin{equation*}
	\Phi_x(t)=\left\{\begin{array}{ll}
		\dfrac{f(t)-f(x)-(t-x)f'(x)-\frac{1}{2}(t-x)^2f''(x)}{(t-x)^2}, & \textrm{ if } t\neq x\\
		0, &\textrm{ if } t=x
	\end{array}\right.
\end{equation*}
Then, it is clear that $\Phi_x(x)=0$. Also, observe that the function $\Phi_x(\cdot)$ belongs to $E$. Hence, by Taylor's theorem, we get 
$$f(t)=f(x)+(t-x)f'(x)+\dfrac{(t-x)^2}{2}f''(x)+(t-x)^2 \Phi_x(t)$$
Now the definition of the operators \eqref{wrightoperators} implies that
$$W_n^{(\beta)}(f,x)-f(x)=f'(x)W_n^{(\beta)}(t-x,x)+\dfrac{1}{2}f''(x)W_n^{(\beta)}((t-x)^2,x)+W_n^{(\beta)}((t-x)^2 \Phi_x(t),x)$$
Therefore, using lemma \ref{moments1}, we have 
\begin{eqnarray}\label{eq4.1}
	\left|n\left(W_n^{(\beta)}(f,x)-f(x)\right)-\frac{1 + 2\beta + 2\beta^2}{2\beta} \, x f''(x)\right| &\leq& n\left|W_n^{(\beta)}((t-x)^2 \Phi_x(t),x)\right|
\end{eqnarray}
Applying the Cauchy-Schwarz inequality to the second term on the right-hand
side of (\ref{eq4.1}), then we see that
$$\left|W_n^{(\beta)}((t-x)^2 \Phi_x(t),x)\right|\leq \sqrt{W_n^{(\beta)}((t-x)^4,x)}\sqrt{W_n^{(\beta)}(\Phi_x^2(t),x)}$$
this yields 
\begin{equation}\label{eq4.2}
	n\left|W_n^{(\beta)}((t-x)^2 \Phi_x(t),x)\right|\leq \sqrt{n^2W_n^{(\beta)}((t-x)^4,x)}\sqrt{W_n^{(\beta)}(\Phi_x^2(t),x)}
\end{equation}
Let $\eta_x(t):=\Phi_x^2(t)$. In this case, observe that $\eta_x(x)=0$ and $\eta_x(\cdot)\in E$. Then it follows from Theorem \ref{theoremstastical} that 
\begin{equation}\label{eq4.3}
	st_A-\lim_{n\to \infty} W_n^{(\beta)}(\Phi_x^2(t),x)=st_A-\lim_{n\to \infty} W_n^{(\beta)}(\eta_x(t),x)=\eta_x(x)=0
\end{equation}
uniformly with respect to $x\in [0,b]$. Now considering \eqref{eq4.2} and \eqref{eq4.3}, and also Lemma \ref{theorem1.23}, we immediately lead to 
\begin{equation}\label{eq4.4}
	st_A-\lim_{n\to \infty} n\left( W_n^{(\beta)}((t-x)^2\Phi_x(t),x)\right)=0
\end{equation}
uniformly with respect to $x\in [0,b]$. Using \eqref{eq4.1} to \eqref{eq4.4} and also considering $\displaystyle st_A-\lim_{n\to \infty} \dfrac{1}{n}=0$, we have 
$$st_A-\lim_{n\to \infty} n\left(W_n^{(\beta)}(f,x)-f(x)\right)\frac{1 + 2\beta + 2\beta^2}{2\beta} \, x f''(x)$$
uniformly with respect to $x \in [0, b]$. 
\end{proof}

\section*{Ethics declarations}
\subsection*{Conflict of interest}
The authors declares that they have no conflict of interest.
\subsection*{Ethical approval}
This article does not contain  any studies with human participants or animals performed by the authors.
\subsection*{Informed consent}
For this type of study informed consent was not required.
\subsection*{Data availability}
Data sharing is not applicable to this article as no new data were created or analysed in this study.


\begin{thebibliography}{10}
	\expandafter\ifx\csname url\endcsname\relax
	\def\url#1{\texttt{#1}}\fi
	\expandafter\ifx\csname urlprefix\endcsname\relax\def\urlprefix{URL }\fi
	\expandafter\ifx\csname href\endcsname\relax
	\def\href#1#2{#2} \def\path#1{#1}\fi
	
	\bibitem{totik1994approximation}
	V.~Totik, Approximation by \textrm{B}ernstein polynomials, American Journal of
	Mathematics 116~(4) (1994) 995--1018.
	
	\bibitem{lewanowicz2004generalized}
	S.~Lewanowicz, P.~Wo{\'z}ny, Generalized \textrm{B}ernstein polynomials, BIT
	Numerical Mathematics 44 (2004) 63--78.
	
	\bibitem{sinha2009inverse}
	T.~Sinha, V.~Gupta, P.~Agrawal, A.~R. Gairola, Inverse theorem for an iterative
	combination of \textrm{B}ernstein-\textrm{D}urrmeyer polynomials, Studia
	Universitatis Babes-Bolyai, Mathematica~(4) (2009).
	
	\bibitem{gupta2013approximation}
	V.~Gupta, Approximation properties by \textrm{B}ernstein--\textrm{D}urrmeyer
	type operators, Complex Analysis and Operator Theory 7 (2013) 363--374.
	
	\bibitem{gupta2021operators}
	V.~Gupta, Operators based on \textrm{P}{\'o}lya distribution and finite
	differences, Mathematical Methods in the Applied Sciences (2021).
	
	\bibitem{mishra2014generalized}
	V.~N. Mishra, P.~Patel, On generalized integral \textrm{B}ernstein operators
	based on $q$-integers, Applied Mathematics and Computation 242 (2014)
	931--944.
	
	\bibitem{ozarslan2013astatistical}
	M.~A. Ozarslan, $a$-statistical convergence of \textrm{M}ittag-\textrm{L}effler
	operators, Miskolc Math. Notes 14~(1) (2013) 209--217.
	
	\bibitem{gupta2020convergence}
	V.~Gupta, Convergence estimates for gamma operator, Bulletin of the Malaysian
	Mathematical Sciences Society 43 (2020) 2065--2075.
	
	\bibitem{gupta2024new}
	V.~Gupta, New operators based on \textrm{L}aguerre polynomials, Revista de la
	Real Academia de Ciencias Exactas, F{\'\i}sicas y Naturales. Serie A.
	Matem{\'a}ticas 118~(1) (2024) 19.
	
	\bibitem{wright1933coefficients}
	E.~M. Wright, On the coefficients of power series having exponential
	singularities, Journal of the London Mathematical Society 1~(1) (1933)
	71--79.
	
	\bibitem{gray1895treatise}
	A.~Gray, E.~Gray, G.~B. Mathews, E.~Meissel, A treatise on \textrm{B}essel
	functions and their applications to physics, Macmillan and Company, 1895.
	
	\bibitem{bateman1953higher}
	H.~Bateman, Higher transcendental functions [volumes i-iii], Vol.~1,
	McGRAW-HILL book company, 1953.
	
	\bibitem{khanmamedov2010zeros}
	A.~K. Khanmamedov, K.~E. Abbasova, On zeros of the modified \textrm{B}essel
	function of the first kind, Azerbaijan Journal of Mathematics 11 (2021)
	2218--6816.
	
	\bibitem{mehrez2017functional}
	K.~Mehrez, Functional inequalities for the \textrm{W}right functions, Integral
	Transforms and Special Functions 28~(2) (2017) 130--144.
	
	\bibitem{mishra2017statistical}
	V.~N. Mishra, P.~Patel, $(\alpha,\beta)$-statistical convergence of modified
	$q$-\textrm{D}urrmeyer operators, Communications Faculty of Sciences
	University of Ankara Series A1 Mathematics and Statistics 66~(2) (2017)
	263--275.
	
	\bibitem{duman2008statistical}
	O.~Duman, A-statistical convergence of sequences of convolution operators,
	Taiwanese Journal of Mathematics 12~(2) (2008) 523--536.
	
	\bibitem{gadjiev2002some}
	A.~Gadjiev, C.~Orhan, Some approximation theorems via statistical convergence,
	The Rocky Mountain Journal of Mathematics (2002) 129--138.
	
	\bibitem{duman2006statistical}
	O.~Duman, E.~Erku{\c{s}}, V.~Gupta, Statistical rates on the multivariate
	approximation theory, Mathematical and Computer Modelling 44~(9-10) (2006)
	763--770.
	
	\bibitem{gupta2009statistical}
	V.~Gupta, C.~Radu, Statistical approximation properties of
	$q$-\textrm{B}askakov-\textrm{K}antorovich operators, Open Mathematics 7~(4)
	(2009) 809--818.
	
	\bibitem{ispir2008statistical}
	N.~Ispir, V.~Gupta, A-statistical approximation by the generalized
	\textrm{K}antorovich-\textrm{B}ernstein type rational operators., Southeast
	Asian Bulletin of Mathematics 32~(1) (2008).
	
	\bibitem{gupta2014convergence}
	V.~Gupta, R.~P. Agarwal, Convergence estimates in approximation theory,
	Vol.~13, Springer, 2014.
	
	\bibitem{gupta2021modifications}
	V.~Gupta, T.~M. Rassias, Modifications of certain operators, in: Computation
	and Approximation, Springer, 2021, pp. 37--68.
	
	\bibitem{gupta2021higher}
	V.~Gupta, Higher order \textrm{L}upa{\c{s}}-\textrm{K}antorovich operators and
	finite differences, Revista de la Real Academia de Ciencias Exactas,
	F{\'\i}sicas y Naturales. Serie A. Matem{\'a}ticas 115~(3) (2021) 1--16.
	
	\bibitem{mursaleen2009invariant}
	M.~Mursaleen, O.~H. Edely, On the invariant mean and statistical convergence,
	Applied Mathematics Letters 22~(11) (2009) 1700--1704.
	
	\bibitem{gurhan2015durrmeyer}
	I.~Gurhan, B.~Cekim, Durrmeyer-type generalization of
	\textrm{M}ittag-\textrm{L}effler operators, Gazi University Journal of
	Science 28~(2) (2015) 259--263.
	
	\bibitem{garrancho2019general}
	P.~Garrancho, A general \textrm{K}orovkin result under generalized convergence,
	Constructive Mathematical Analysis 2~(2) (2019) 81--88.
	
	\bibitem{chandra2022approximation}
	P.~Chandra, V.~D. Kumar, D.~Naokant, Approximation by durrmeyer variant of
	\textrm{C}heney-\textrm{S}harma \textrm{C}hlodovsky operators, Mathematical
	Foundations of Computing (2022) 0--0\href
	{https://doi.org/10.3934/mfc.2022034} {\path{doi:10.3934/mfc.2022034}}.
	
	\bibitem{agrawal2022stancu}
	P.~N. Agrawal, S.~Singh, Stancu variant of
	\textrm{J}akimovski-\textrm{L}eviatan-\textrm{D}urrmeyer operators involving
	\textrm{B}renke type polynomials, Mathematical Foundations of Computing
	(2022).
	
	\bibitem{maurya2018approximation}
	R.~Maurya, H.~Sharma, C.~Gupta, Approximation properties of
	\textrm{K}antorovich type modifications of $(p, q)-$
	\textrm{M}eyer-\textrm{K}{\"o}nig-\textrm{Z}eller operators, Constructive
	Mathematical Analysis 1~(1) (2018) 58--72.
	
	\bibitem{duman2004statistical}
	O.~Duman, C.~Orhan, Statistical approximation by positive linear operators,
	Studia Mathematica 161~(2) (2004).
	
	\bibitem{karaisa2017alphabeta}
	A.~Karaisa, U.~Kadak, On $\alpha$$\beta$-statistical convergence for sequences
	of fuzzy mappings and \textrm{K}orovkin type approximation theorem, Filomat
	31~(12) (2017) 3749--3760.
	
	\bibitem{bozma2022approximation}
	G.~Bozma, N.~G. Bilgin, Approximation in a variable bounded interval, in: 8th
	International Conference on Intuitionistic Fuzzy Sets and Contemporary
	Mathematics, 2022, p.~25.
	
	\bibitem{karaisaa2017alphabeta}
	A.~Karaisaa, U.~Kadakb, On $\alpha$$\beta$-statistical convergence for
	sequences of fuzzy mappings and \textrm{K}orovkin type approximation theorem,
	Filomat 31~(12) (2017) 3749--3760.
	
	\bibitem{gogoi2021weighted}
	J.~Gogoi, H.~Dutta, Weighted ($\lambda$, $\mu$)-statistical convergence and
	statistical summability methods of double sequences of fuzzy numbers with
	application to \textrm{K}orovkin-type fuzzy approximation theorem, Soft
	Computing 25~(11) (2021) 7645--7656.
	
	\bibitem{karakaya2015korovkin}
	V.~Karakaya, A.~Karaisa, Korovkin type approximation theorems for weighted
	$\alpha\beta$-statistical convergence, Bulletin of Mathematical Sciences
	5~(2) (2015) 159--169.
	
	\bibitem{braha2020korovkin}
	N.~L. Braha, V.~Loku, Korovkin type theorems and its applications via
	$\alpha\beta$-statistically convergence, J. Math. Inequal 14~(4) (2020)
	951--966.
	
	\bibitem{abay2018alphabeta}
	M.~Abay~Tok, E.~E. Kara, S.~Altundag, On the $\alpha\beta$-statistical
	convergence of the modified discrete operators, Advances in Difference
	Equations 2018~(1) (2018) 1--6.
	
	\bibitem{aktuuglu2016weighted}
	M.~A. {\''O}zarslan, H.~Aktu{\u{g}}lu, Weighted $\alpha$$\beta$-statistical
	convergence of \textrm{K}antorovich-\textrm{M}ittag-\textrm{L}effler
	operators, Mathematica Slovaca 66~(3) (2016) 695--706.
	
	\bibitem{fast1951convergence}
	H.~Fast, Sur la convergence statistique, in: Colloquium mathematicae, Vol.~2,
	1951, pp. 241--244.
	
	\bibitem{aktuuglu2014korovkin}
	H.~Aktu{\u{g}}lu, Korovkin type approximation theorems proved via
	$\alpha$$\beta$-statistical convergence, Journal of Computational and Applied
	Mathematics 259 (2014) 174--181.
	
	\bibitem{patel2024Wright}
	P.~G. Patel, On positive linear operators linking gamma,
	\textrm{M}ittag-\textrm{L}effler and \textrm{W}right functions, submitted to
	jounral.
	
	\bibitem{patel2024Wright1}
	P.~G. Patel, On \textrm{W}right operators involving generalized gamma function,
	submitted to jounral.
	
\end{thebibliography}
\end{document}